\documentclass{article}

\title{RATIONAL PERIOD FUNCTIONS ON THE HECKE GROUPS
		\thanks{1991 Mathematics Subject Classification. 
				Primary 11F67;
				Secondary 11F12, 11E45.
				}
	  }

\author{Wendell Culp-Ressler \\
	  Department of Mathematics \\
	  Franklin \& Marshall College \\
	  Lancaster, PA 17604-3003}

\date{July 11, 2000}

\usepackage{amsmath,amsthm,amsfonts,amssymb}            
\usepackage{makeidx}

\makeindex

   \theoremstyle{plain} 				
   \newtheorem{theorem}{Theorem}
   \newtheorem*{corollary}{Corollary}
   \newtheorem{lemma}{Lemma}

   \theoremstyle{definition}
   \newtheorem{defn}{Definition}

   \theoremstyle{remark}
   \newtheorem*{remark}{Remark}

   \newcommand{\Z}{\mathbb{Z}}  	
   \newcommand{\Q}{\mathbb{Q}}  	
   \newcommand{\R}{\mathbb{R}}  	
   \newcommand{\C}{\mathbb{C}}  	
   
   \newcommand{\uhp}{\mathcal{H}}	
   \newcommand{\Pset}{P}        	

\begin{document}

\maketitle

\begin{abstract}
We describe rational period functions on the Hecke 
groups and characterize the ones whose poles satisfy a certain symmetry.
This generalizes part of the characterization of rational period 
functions on the modular group, which is one of the Hecke groups.
\end{abstract}

\section{INTRODUCTION}

Marvin Knopp \cite{Kno78} introduced the idea of a rational period 
function (RPF) for an automorphic integral of weight \(2k\) \((k\in \Z)\) on 
\(\Gamma\), where \(\Gamma\) is any Fuchsian group acting on the upper 
half-plane.

Knopp \cite{Kno78,Kno81}, Hawkins \cite{Haw}, and Choie and Parson 
\cite{CP90,CP91} took steps toward an explicit characterization of RPFs on 
the modular group \(\Gamma(1)\).
Ash \cite{Ash89} used cohomological techniques to provide a 
characterization, after which Choie and Zagier \cite{CZ93} and, 
independently, Parson \cite{Par93} provided an explicit characterization 
of the RPFs on \(\Gamma(1)\).
The explicit characterizations use negative continued fractions 
to establish a connection between the poles of 
RPFs and binary quadratic forms.

Schmidt \cite{Sch96} generalized Ash's work, giving an abstract 
characterization of RPFs on any finitely generated Fuchsian group of the 
first kind with parabolic elements, a class of groups which includes the 
Hecke groups.
Schmidt \cite{Sch93} and Schmidt and Sheingorn \cite{SS95} have taken 
steps toward an explicit 
characterization of RPFs on the Hecke groups using generalizations of the 
classical continued fractions and binary quadratic forms.

In this paper we give an explicit characterization of a class of 
rational period functions on the Hecke groups.
Our characterization is for rational period functions of weight 
\( 2k \), with \( k \) a positive odd integer, 
and with irreducible pole 
sets which are Hecke-symmetric, as defined in Section \ref{subsec.BQFs}.
We also determine an explicit expression for all rational period 
functions on Hecke groups.
We use the correspondence, developed in \cite{C-R99a},
between poles of rational 
period functions and generalized binary quadratic forms.

\section{DEFINITIONS AND BACKGROUND IDEAS}
\label{sec.BackIdeas}

\subsection{Hecke groups}
\label{subsec.Heckegroups}

Let \( \lambda \) be a positive real number and
let \(G(\lambda)=\langle S,T \rangle / \{\pm I\}\),
where 
\(S = S_{\lambda} 
   = \bigl( \begin{smallmatrix} 1&\lambda\\ 0&1 \end{smallmatrix} \bigr)\),
\(T = \bigl( \begin{smallmatrix} 0&-1\\ 1&0 \end{smallmatrix} \bigr)\),
and \(I\) is the identity.
Erich Hecke \cite{Hec36} showed that the only 
values of \(\lambda\) between \(0\) and \(2\) for which \(G(\lambda)\) is discrete 
are 
\(\lambda = \lambda_{p}
         = 2\cos(\pi/p)\) 
for \(p=3, 4, 5, \dotsc\).
These are the 
    \index{Hecke group}%
\emph{Hecke groups} which we will denote by
\(G_{p}=G(\lambda_{p})\) for \(p\geq3\).
Put \( U = U_{\lambda_{p}} = S_{\lambda_{p}}T \in G_{p} \).
The Hecke groups have two relations, 
\(T^{2}=U^{p}=I\).

Hecke showed that \(G(\lambda)\) is also discrete if \(\lambda \geq 2\),
however these groups have only one relation \(T^{2}=I\).
In this case the related rational period functions have a simpler 
structure, which is given in \cite{HK92}.

For 
\(M = \bigl( \begin{smallmatrix} a&b\\ c&d \end{smallmatrix} \bigr)
    \in G_{p}\),
we have \(a, b, c, d \in \Z [\lambda_{p}]\) and \(ad-bc=1\), 
so \(G_{p}\) is a subgroup of \(\textrm{SL}(2,\Z [\lambda_{p}])\).
It is well-known that \(G_{3}=\textrm{SL}(2,\Z[\lambda_{3}])\) 
(\emph{i.e.},\(\Gamma(1)=\textrm{SL}(2,\Z)\)),
however for the other Hecke groups 
\(G_{p} \subsetneqq \textrm{SL}(2,\Z[\lambda_{p}])\).

Members of Hecke groups 
    \index{Hecke group!action on \( \C \)}%
act on the Riemann sphere as linear fractional transformations.
An element 
\( M = \bigl( \begin{smallmatrix} a&b\\ c&d \end{smallmatrix} \bigr)
      \in G_{p} \) 
is
    \index{hyperbolic!matrix}%
\emph{hyperbolic} if \( \vert a+d \vert >2 \),
    \index{parabolic matrix}%
\emph{parabolic} if \( \vert a+d \vert =2 \),
and 
    \index{elliptic matrix}%
\emph{elliptic} if \( \vert a+d \vert <2 \). 
We will designate fixed points accordingly.

We will use the following lemma 
when we study the poles of rational 
period functions.
\begin{lemma}\label{lem.PositiveEntries}
  Fix \( p \geq 3 \) and let \( U = U_{\lambda_{p}} \).
  The nonzero entries of \( U^{n}T \in G_{p}\) are 
  all positive for \( 1 \leq n \leq p-1 \).  
  The only zero entries occur in 
  \( UT = 
    \bigl( \begin{smallmatrix} 
      1&\lambda\\ 0&1 
    \end{smallmatrix} \bigr) 
  \) 
  and
  \( U^{p-1}T = 
    \bigl( \begin{smallmatrix} 
      1&0\\ \lambda&1 
    \end{smallmatrix} \bigr).
  \)
\end{lemma}

\begin{proof}
    Write 
    \( U^{n} = 
      \bigl( \begin{smallmatrix} 
        a_{n}&b_{n}\\ c_{n}&d_{n} 
      \end{smallmatrix} \bigr).
    \)
    By induction on \( n \) we have that 
    \( a_{n} = \frac{sin\left((n+1)\pi/p\right)}{sin(\pi/p)} \)
    for \( n \geq 0 \), 
    from which it is clear that \( a_{n}>0 \) for 
    \( 0 \leq n \leq p-2 \).
    We write \( U^{n} = UU^{n-1} = U^{n-1}U \), or
    \[ U^{n} =   
      \begin{pmatrix} 
        \lambda a_{n-1}-c_{n-1}&\lambda b_{n-1}-d_{n-1}\\ a_{n-1}&b_{n-1} 
      \end{pmatrix} 
      =
      \begin{pmatrix} 
        \lambda a_{n-1}+b_{n-1}& -a_{n-1}\\ \lambda c_{n-1}+d_{n-1}&-c_{n-1} 
      \end{pmatrix}.
   \]
   Then 
   \begin{equation*}
       U^{n} =   
          \begin{pmatrix} 
            a_{n}& -a_{n-1}\\ a_{n-1}&-a_{n-2} 
          \end{pmatrix},
   \end{equation*}
   so
   \[ U^{n}T = 
      \begin{pmatrix} 
        a_{n-1}& a_{n}\\ a_{n-2}&a_{n-1} 
      \end{pmatrix},
   \]
   which has positive entries for \( 2 \leq n \leq p-2 \).
\end{proof}

\subsection{Rational period functions}
\label{subsec.RPFs}

For 
\(M = \bigl( \begin{smallmatrix} *&*\\ c&d \end{smallmatrix} \bigr) 
\in \textrm{SL}(2,\R)\)
and \(f(z)\) a complex function,
we define the 
    \index{weight}%
\emph{weight \( 2k \) 
    \index{slash operator}%
slash operator}
\(f \mid_{2k} M = f \mid M\) by
\begin{equation*}
	\left( f\mid M \right) (z) = (cz+d)^{-2k} f \left( Mz \right).
\end{equation*}
\begin{defn}
Let \( U = U_{\lambda_{p}}\) for \(p\geq 3\) and let \(k \in \Z^{+} \).
A 
    \index{rational period function}%
\emph{rational period function} (RPF) 
\emph{of     
    \index{weight}%
weight \(2k\) for 
\(G_{p}\)} is a 
rational function \(q(z)\) which satisfies the relations
\begin{equation}
	q+q\mid T = 0,
  \label{eq:FirstRelation}
\end{equation}
and 
\begin{equation}
	q+q\mid U + q\mid U^{2} + \dots + q\mid U^{p-1} = 0.
\label{eq:SecondRelation}
\end{equation}
\label{def:RPFbyRelations}
\end{defn}

Marvin Knopp \cite{Kno78} first introduced RPFs as period functions for 
automorphic integrals in the following way.  
\begin{defn}\label{def:RPFbyMI}
    Let \( S = S_{\lambda_{p}}\) for \(p\geq 3\), and let 
    \(k \in \Z^{+} \).
    Suppose that \(F\) is a function meromorphic in \(\uhp\) and at \(i\infty\) 
    which satisfies
    \begin{equation*}
	\left( F\mid S \right) (z) = F(z),
    \end{equation*}
    and 
    \begin{equation*}
	\left( F\mid T \right) (z) = F(z) + q(z),
    \end{equation*}
    where \(q(z)\) is a rational function.
    Then \(F\) is an
        \index{automorphic integral}%
    \emph{automorphic integral of 
        \index{weight}%
    weight \(2k\) on \( G_{p} \) 
    with 
        \index{rational period function}%
    rational period function \(q(z)\)}.
\end{defn}
These definitions are equivalent.
It is easy to show that any rational period function satisfying 
Definition \ref{def:RPFbyMI} also satisfies Definition \ref{def:RPFbyRelations}.
On the other hand, Knopp showed \cite[Theorem 3]{Kno74} that any function 
satisfying Definition \ref{def:RPFbyRelations} is the period function for an 
automorphic integral of weight \(2k\).

For any rational period function \(q\) on \(G_{p}\) 
we let 
    \index{\(P(q) \)}%
    \index{poles!\(P(q) \)}%
\( P(q) \) denote the set of poles of \(q\).
Hawkins \cite{Haw} introduced the idea of an 
    \index{irreducible system of poles}%
\emph{irreducible system of poles}
(or
    \index{irreducible pole set}%
\emph{irreducible pole set}),
the minimal set of poles which are forced to occur together by the 
relations \eqref{eq:FirstRelation} and \eqref{eq:SecondRelation}.
If \( \alpha \in \Pset^{*}(q) = \Pset(q) \setminus \{0\} \), 
we let \( \Pset(q;\alpha) \) 
    \index{\(P(q;\alpha) \)}%
    \index{poles!\(P(q;\alpha) \)}%
denote the irreducible set of poles of 
\( q \) which contains \( \alpha \).

\subsection{Binary quadratic forms}
\label{subsec.BQFs}

An
    \index{indefinite binary quadratic form}%
    \index{binary quadratic form}%
\emph{indefinite binary quadratic form} on \( \Z[\lambda_{p}] \)
is an expression of the form
\[ Q(x,y) = Ax^{2}+Bxy+Cy^{2}, \]
with \( A, B, C \in \Z[\lambda_{p}] \)
and \( D=B^{2}-4AC > 0 \).
We denote such a form by \( Q=[A,B,C] \) and refer to it as a
\( \lambda_{p} \)-BQF.
If \( \Z[\lambda_{p}] \) is a principal ideal domain we also require 
that a \( \lambda_{p} \)-BQF \( [A,B,C] \) be 
    \index{primitive!BQF}%
\emph{primitive}, \emph{i.e.}, that 
\( (A) + (B) + (C) = (1) \), where 
\( (x) = x\Z[\lambda_{p}] \) denotes the ideal of \( \Z[\lambda_{p}] \) 
generated by \( x \).

We need a \( 1-1 \) correspondence between \( \lambda_{p} \)-BQFs and 
certain algebraic numbers.
We map the \( \lambda_{p} \)-BQF \( Q = [A,B,C] \) to the  number
\( \alpha_{Q} = \frac{-B+\sqrt{D}}{2A} \in \Q(\lambda_{p},\sqrt{D}) \).
    \index{\(\alpha_{Q} \)}%
In this form the mapping is invertible only 
for \( p = 3 \),  \( G_{p} = \Gamma(1) \).
When \( p > 3 \),
ambiguity arises from the presence of nontrivial units in 
\( \Z[\lambda_{p}] \)
and from the fact that the discriminant \( D \) 
need not be square-free in \( \Z[\lambda_{p}] \).
In order to recover a \( 1-1 \) correspondence from this mapping we 
first restrict the range of the 
mapping to hyperbolic fixed points of \( G_{p} \).
Then we give an algorithm which produces a unique inverse for any 
hyperbolic number.

\begin{lemma}\label{lem.NumbertoForm}
    If \( \alpha \) is a hyperbolic fixed point of \( G_{p} \) 
    it may be associated with a unique indefinite 
    \( \lambda_{p} \)-BQF
    \( Q = Q_{\alpha} \)
        \index{\(Q_{\alpha} \)}%
    such that
    \( \alpha = \alpha_{Q} \).
    If \( \Z[\lambda_{p}] \) is a principal ideal domain,
    \( Q \) may be chosen to be primitive.
\end{lemma}
If an indefinite \( \lambda_{p} \)-BQF \( Q \) is associated 
with a hyperbolic 
number as in the Lemma, 
we say that \( Q \) is
    \index{hyperbolic!BQF}%
\emph{hyperbolic}.
\begin{proof}
    We present an outline of the more detailed proof in \cite{C-R99a}.
    Suppose that 
    \( \alpha \) is a hyperbolic fixed point of \( G_{p} \).
    Let
    \(  M_{\alpha} = \bigl( \begin{smallmatrix} a&b\\ c&d \end{smallmatrix} \bigr)
        \in G_{p} \) 
    be a generator of the cyclic group of matrices fixing
    \( \alpha \).
    \( M_{\alpha} \) is determined up to 
    inverses and
    \[ \alpha = \frac{a-d \pm \sqrt{(a+d)^{2}-4}}{2c}. \]
    We show in \cite{C-R99a} that if \( \Z[\lambda_{p}] \) is a principal 
    ideal domain, there is a uniquely determined positive number
    \( g \in \Z[\lambda_{p}] \) with \( (c,d-a,-b) = (g) \).
    Then \( \frac{1}{g}[c,d-a,-b] \) is primitive and the unique 
    primitive \( \lambda_{p} \)-BQF for \( \alpha \) is 
    \begin{equation*}
	Q_{\alpha} = 
	\begin{cases}
	    \frac{1}{g}[c,d-a,-b], & 
	        \text{if } \alpha = \frac{a-d+\sqrt{(a+d)^{2}-4}}{2c}, \\
	    \frac{-1}{g}[c,d-a,-b], & 
	        \text{if } \alpha = \frac{a-d-\sqrt{(a+d)^{2}-4}}{2c}.
	\end{cases}
    \end{equation*}
    If \( \Z[\lambda_{p}] \) fails to be a principal ideal domain
    we put \( g=1 \).
\end{proof}

Suppose that \( \alpha = \frac{-B+\sqrt{D}}{2A} \) is a hyperbolic point 
associated with the \( \lambda_{p} \)-BQF
\( Q_{\alpha} = [A,B,C] \).
We define the 
    \index{Hecke!conjugate}%
\emph{Hecke} (\( \lambda_{p} \)-)\emph{conjugate of} \( \alpha \) 
to be
\( \alpha^{\prime} = \frac{-B-\sqrt{D}}{2A} \),
the number associated with 
\( -Q_{\alpha} = [-A,-B,-C] \),
\emph{i.e.}, \( \alpha_{Q}^{\prime}=\alpha_{-Q} \).
A calculation shows that if \( V \in G_{p} \),
then \( \left(V\alpha\right)^{\prime} = V\alpha^{\prime} \).
In the case of \( G_{3}=\Gamma(1) \), Hecke conjugation reduces to 
algebraic conjugation over \( \Q \).

If \( R \) is a set of hyperbolic fixed points of 
\( G_{p} \) we write
\( R^{\prime} = \{x^{\prime} \mid x \in R\} \).
We say that \( R \) has 
    \index{Hecke!symmetry}%
\emph{Hecke} (\( \lambda_{p} \)-)\emph{symmetry}
if \( R = R^{\prime} \).

Elements of a Hecke group
    \index{Hecke group!action on BQFs}%
act on \( \lambda_{p} \)-BQFs by
\( \left(Q \circ M \right)(x,y) = Q(a x + b y, c x + d y) \)
for 
\( M = \bigl( \begin{smallmatrix} a &b \\ c & d 
              \end{smallmatrix} \bigr)\ \in G_{p} \).
This action preserves the discriminant and
maps primitive forms to primitive forms.
We say that \( \lambda_{p} \)-BQFs \( Q_{1} \) and \( Q_{2} \) are
    \index{equivalent!BQFs}%
\emph{\( G_{p} \)-equivalent},
and write \( Q_{1} \sim Q_{2} \),
if there exists a \( V \in G_{p} \) such that 
\( Q_{2} = Q_{1} \circ V \).
It is easy to check that \( G_{p} \)-equivalence is an equivalence relation, so 
\( G_{p} \) partitions the \( \lambda_{p} \)-BQFs into equivalence 
classes of forms.

Every \( G_{p} \)-equivalence class of \( \lambda_{p} \)-BQFs contains 
either all hyperbolic forms or no hyperbolic forms \cite{C-R99a}.
Consequently, we may label equivalence classes themselves as hyperbolic or 
non-hyperbolic.

Suppose that \( Q_{1} \) and \( Q_{2} \) are hyperbolic 
\( \lambda_{p} \)-BQFs.
A calculation shows that \( Q_{2} = Q_{1} \circ V \) 
if and only if \( \alpha_{2} = V^{-1}\alpha_{1} \), 
where \( \alpha_{1} \) and \( \alpha_{2} \) are associated with 
\( Q_{1} \) and \( Q_{2} \), respectively.
Thus \( G_{p} \)-equivalence of hyperbolic \( \lambda_{p} \)-BQFs 
induces a corresponding 
    \index{equivalent!numbers}
\( G_{p} \)-equivalence of associated numbers.

If \( \mathcal{A} \) denotes a \( G_{p} \)-equivalence 
class of \( \lambda_{p} \)-BQFs, we let
    \index{\(-\mathcal{A}\)}
\( -\mathcal{A} 
     = \left\{-Q \vert Q \in \mathcal{A}\right\} \).
Then \( -\mathcal{A} \) is another \( G_{p} \)-equivalence 
class of forms, not necessarily distinct from \( \mathcal{A} \).
If \( \mathcal{A} \) is hyperbolic, so is \( -\mathcal{A} \), and 
the numbers associated with the forms in 
\( -\mathcal{A} \) are the Hecke 
conjugates of the numbers associated with the forms in 
\( \mathcal{A} \).

We say that a hyperbolic \( \lambda_{p} \)-BQF
\( Q=[A,B,C] \) is 
    \index{simple!BQF}
(\( G_{p} \)-)\emph{simple} if \( A>0>C \).
If \( Q \) is \( G_{p} \)-simple  
we say that the associated hyperbolic number
\( x_{Q} \) is 
    \index{simple!number}
\mbox{(\( G_{p} \)-)\emph{simple}}.
We show in \cite{C-R99a} that a hyperbolic fixed point \( x \) is 
\( G_{p} \)-simple if and only if \( x^{\prime} < 0 < x \).

\section{POLES OF RATIONAL PERIOD FUNCTIONS}
\label{sec.RPF_Poles}

Throughout this section we fix \( p \geq 3 \), 
\( \lambda = \lambda_{p} \), and
\( U = U_{\lambda_{p}} \).
We assume that \(q\) an RPF of weight \(2k \in 2\Z^{+} \) on \(G_{p}\)
with pole set \(\Pset = \Pset(q) \).

\subsection{Sets of Poles}
\label{subsec.PoleSets}

In this section we show that the positive poles of 
\( q \) may be put into cycles, 
which we use to write the irreducible systems of poles.
We also make the connection between these cycles of poles and simple
\( \lambda_{p} \)-BQFs.

\begin{lemma}\label{lem.PolesReal}
  \(\Pset(q) \subset \R\) for any RPF \( q \) of weight 
  \( 2k \in 2 \Z^{+} \) on \( G_{p} \).
\end{lemma}

\begin{remark}
  This is given without proof by Meier and Rosenberger \cite{MR84},
  who state that the proof is analogous to the one by Knopp \cite{Kno81} for 
  rational period functions on the modular group.
  The following proof is a generalization of the one by Choie and Zagier
  \cite{CZ93} for the modular group.
\end{remark}

\begin{proof}
Suppose, by way of contradiction, that \(\alpha_{1} \in \Pset \) but
\(\alpha_{1} \notin \R\).
By the first relation \eqref{eq:FirstRelation} we have 
\(T\alpha_{1} \in \Pset \), 
then by the second relation \eqref{eq:SecondRelation} we have
\(U^{j_{1}}T\alpha_{1} \in \Pset \) for some 
\(j_{1}\), \(1\leq j_{1}\leq p-1\).
Repeating this, we produce a sequence of poles 
\(\{\alpha_{1}, \alpha_{2}, \dots \}\), none of which are real, with 
\(\alpha_{\nu +1}=U^{j_{\nu}}T\alpha_{\nu}\) and \(1\leq j_{\nu} \leq p-1\)
for each \(\nu \geq 1\).
Lemma \ref{lem.PositiveEntries}  and
a geometric argument show that 
\(\left| \arg(\alpha_{\nu}) \right| > \left| \arg(\alpha_{\nu +1}) \right|\) 
for each \(\nu \geq 1\).
But this is impossible, since \(\Pset\) is a finite set.
\end{proof}

Put \( \lambda = \lambda_{p} \) and \( U = U_{\lambda} \) 
for \( p \geq 3 \).
Define 
\( \Phi_{p}: [0,\infty) \rightarrow [0,\infty) \)
    \index{\(\Phi_{p}\)}%
by
\[ 
  \Phi_{p}(x) = 
      \begin{cases}
          TUx, & U^{p}(0) \leq x < U^{p-1}(0) \\
          TU^{2}x, & U^{p-1}(0) \leq x < U^{p-2}(0) \\
	  \vdots \\
	  TU^{p-1}x,   & U^{2}(0) \leq x.
      \end{cases}
    \]
This function is given more explicitly by
\[ 
  \Phi_{p}(x) = 
      \begin{cases}
          \frac{x}{1-\lambda x},   
	              & 0 \leq x < 1/\lambda \\
	  \frac{1-\lambda x}{(\lambda^{2}-1)x-\lambda },   
	              & 1/\lambda \leq x < \lambda/(\lambda^{2}-1) \\
	  \vdots \\
	  x-\lambda,   & \lambda \leq x.
      \end{cases}
    \]
It can be shown that \( TU^{n}(x) > 0 \) if and only if
\( U^{p-n+1}(0) < x < U^{p-n}(0) \),
so the exponent \( n \) in 
\( \Phi_{p}(x) = TU^{n}x \) is the \emph{unique} exponent 
between \( 1 \) and \( p-1 \) for which \( TU^{n}x > 0 \)
\cite{C-R99a}.

If \( \mathcal{A} \) is a hyperbolic \( G_{p} \)-equivalence class of 
\( \lambda_{p} \)-BQFs we write
    \index{\(Z_{\mathcal{A}}\)}%
\linebreak[1]
\( Z_{\mathcal{A}} 
    = \left\{x \mid Q_{x} \in \mathcal{A}, Q_{x} \textrm{ simple} \right\} \).
We show in \cite{C-R99a} that each hyperbolic equivalence class 
of \( \lambda_{p} \)-BQFs contains at 
least one simple form, 
\emph{i.e.}, 
\( Z_{\mathcal{A}} \neq \emptyset \) 
for every hyperbolic \( \mathcal{A} \).
We also show in \cite[Theorem 3]{C-R99a} that
the finite cyclic orbits of 
\( \Phi_{p} \) 
are the sets \( \{0\} \)
and \( Z_{\mathcal{A}} \), where 
\( \mathcal{A} \) runs over all hyperbolic 
\( G_{p} \)-equivalence classes of 
\( \lambda_{p} \)-BQFs.

Let 
    \index{\(P^{+} \)}%
    \index{poles!\(P^{+} \)}%
\( \Pset^{+} \) denote the set of positive poles in \( \Pset \),
    \index{\(P^{-} \)}%
    \index{poles!\(P^{-} \)}%
\( \Pset^{-} \) denote the set of negative poles in \( \Pset \),
so \( \Pset^{*} = \Pset^{+} \cup \Pset^{-} \).
Write 
    \index{\(T\mathcal{Z_{\mathcal{A}}} \)}
\( T\mathcal{Z_{\mathcal{A}}} 
    = \left\{T\alpha \mid \alpha \in \mathcal{Z_{\mathcal{A}}} \right\} \).

The following Lemma establishes the relationship between 
irreducible pole sets of RPFs on 
\( G_{p} \) and equivalence classes of simple
\( \lambda_{p} \)-BQFs.
\begin{lemma}\label{lem.MinPoleSet}
    Suppose that \( \alpha \) is a positive pole of an RPF on \( G_{p} \).
    Then 
    \( \alpha \) is hyperbolic and
    \( \alpha \in Z_{\mathcal{A}} \),
    where \( \mathcal{A} \) is the 
    \( G_{p} \)-equivalence class of \( \lambda_{p} \)-BQFs containing
    \( Q_{\alpha} \).
    The irreducible pole set containing \( \alpha \) is 
    \begin{equation*}
        P(q;\alpha) 
	    = Z_{\mathcal{A}} \cup TZ_{\mathcal{A}}.
    \end{equation*}
\end{lemma}

\begin{proof}
    Suppose that \(\alpha = \alpha_{1} \in \Pset^{+}\).
    As in the proof of Lemma \ref{lem.PolesReal}, 
    \(\alpha_{2}=U^{j_{1}}T\alpha_{1} \in \Pset\) for some \(j_{1}\), 
    \(1 \leq j_{1} \leq p-1\).
    By Lemma \ref{lem.PositiveEntries} we may take each
    entry of \(U^{j_{1}}T\) to be non-negative,
    which implies that 
    \(\alpha_{2}=U^{j_{1}}T\alpha_{1}>0\).
    Repeating this process gives a sequence of poles 
    \(\{ \alpha_{1}, \alpha_{2}, \dots \} \subseteq P^{+} \), with 
    \begin{equation}
	\alpha_{\nu +1} = U^{j_{\nu}}T\alpha_{\nu},
	\label{eq:NextPole}
    \end{equation}
    and \(1\leq j_{\nu}\leq p-1\)
    for each \(\nu \geq 1\).
    Since \(\Pset^{+}\) is finite we must have that 
    \begin{equation}
	\alpha_{1} = \alpha_{r+1} = U^{j_{r}}T\alpha_{r},
	\label{eq:FirstPole}
    \end{equation}
    for some \( r \geq 1 \).
    Thus we have a finite cycle of positive poles,
    \( \{\alpha_{1}, \alpha_{2}, \dots, \alpha_{r} \} \).
    Reversing \eqref{eq:NextPole} we have that 
    \( \alpha_{r} = TU^{p-j_{r}}\alpha_{1} 
        = \Phi_{p}(\alpha_{\nu + 1}) \),
    since \( 1 \leq p-j_{\nu} \leq p-1 \) and \( \alpha_{\nu} > 0 \)
    for each \( \nu \geq 1 \).
    Reversing \eqref{eq:FirstPole} we have that 
    \( \alpha_{\nu} = TU^{p-j_{\nu}}\alpha_{\nu + 1}  
        = \Phi_{p}(\alpha_{1}) \),
    since \( 1 \leq p-j_{r} \leq p-1 \) and \( \alpha_{r} > 0 \).
    Hence 
    \( \{\alpha_{1}, \alpha_{2}, \dots, \alpha_{r} \} \subseteq P^{+} \)
    is a finite cyclic orbit of \( \Phi_{p} \), so 
    \( \{\alpha_{1}, \alpha_{2}, \dots, \alpha_{r} \} = Z_{\mathcal{A}} \), 
    where \( \mathcal{A} \) is some
    \( G_{p} \)-equivalence class of 
    \( \lambda_{p} \)-BQFs.
    Hence \( \mathcal{A} \) is hyperbolic and contains
    \( Q_{\alpha} \),
    while \( Q_{\alpha} \) and \( \alpha \) 
    are both hyperbolic and simple.

    Since \( \alpha \) is a positive pole,
    each element of
    \( Z_{\mathcal{A}} \) must also be positive pole, so
    \( Z_{\mathcal{A}} \subseteq P(q;\alpha) \).
    This, with the first relation \eqref{eq:FirstRelation},
    implies that every element of \( TZ_{\mathcal{A}} \) 
    is a negative pole,
    so \( TZ_{\mathcal{A}} \subseteq P(q;\alpha) \).
    On the other hand, if
    \( \beta \in P(q;\alpha) \),
    then \( \beta = M\alpha \),
    where \( M \in G_{p} \) is a product of matrices, each one equal to 
    \( T \) or to \( U^{i} \), \( 1 \leq i \leq p-1 \).
    Thus \( \alpha \) and \( \beta \) are \( G_{p} \)-equivalent 
    numbers, hence 
    \( Q_{\alpha} \) and \( Q_{\beta} \) are \( G_{p} \)-equivalent 
    BQFs,
    so \( Q_{\beta} \in \mathcal{A} \).
    If \( \beta > 0 \), then  
    \( \beta \in Z_{\mathcal{A}} \).
    If \( \beta < 0 \), then 
    \( T\beta \in Z_{\mathcal{A}} \),
    and \( \beta \in TZ_{\mathcal{A}} \).
\end{proof}

\begin{corollary}
    The set of nonzero poles of an RPF on \( G_{p} \) 
    has the form
    \begin{equation}
        \label{eq:PoleDecomp2}
	P^{*} = \bigcup_{\ell=1}^{L}
	        \left( Z_{\mathcal{A}_{\ell}} 
                    \cup TZ_{\mathcal{A}_{\ell}} \right)
    \end{equation}
    where the \( \mathcal{A}_{\ell} \) are distinct indexed hyperbolic 
    \( G_{p} \)-equivalence classes of \( \lambda_{p} \)-BQFs.
\end{corollary}

\begin{remark}
    This is essentially Lemma 3.2 in \cite{Sch93}.
\end{remark}

The following Lemma allows us to express the negative poles 
in \( \Pset^{*} \)
using Hecke conjugation instead of the action of \( T \).
\begin{lemma}\label{lem.NegPoles}
    Let \( \mathcal{A} \) be a \( G_{p} \)-equivalence class of 
    \( \lambda_{p} \)-BQFs.
    Then 
    \( TZ_{\mathcal{A}} = 
    Z_{-\mathcal{A}}^{\prime} \).
\end{lemma}
\begin{proof}
    A routine exercise showing containment in both directions 
    establishes the Lemma.
\end{proof}

\subsection{Principal Parts}
\label{subsec.PPs}

In this section we determine the principle part at any pole of 
any RPF of weight \( 2k \) on \( G_{p} \).
We use this to give an explicit expression of the form of any RPF.

\subsubsection{The Pole at Zero}
\label{subsubsec.PoleatZero}

\begin{lemma}\label{lem.PoleatZero}
    Suppose that \( q \) is an RPF of weight \( 2k \) on \( G_{p} \).
    Then
    \begin{enumerate}
    \renewcommand{\labelenumi}{(\roman{enumi})}
        \item  \( q \) is regular at \( \infty \),
    
        \item  if \( q \) has a pole at \( 0 \), the order of the pole 
        is at most \( 2k \), and
    
        \item  \( q \) has a pole at \( 0 \) of order \( 2k \) 
	    if and only if 
	    \( q(\infty) \neq 0 \).
    \end{enumerate}
\end{lemma}
\begin{proof}
    Fix \( p \geq 3 \) and put \( \lambda=\lambda_{p} \) and \( U = U_{\lambda_{p}} \).
    For \emph{(i)} we apply \( \vert S \) to the second relation
    \eqref{eq:SecondRelation}, use the first relation
    \eqref{eq:FirstRelation} and rearrange to get
    \begin{equation*}
        \left(q-q\vert S\right)(z) 
	    = \left(q\vert US + q\vert U^{2}S + \cdots  + q\vert U^{p-2}S
	      \right)(z).
    \end{equation*}
    We claim that the right hand side of this expression
    approaches \( 0 \) as \( z \rightarrow \infty \).
    To prove our claim it suffices to show that for 
    \( 1 \leq n \leq p-2 \),
    \begin{displaymath}
        \left( q \vert U^{n}S \right)(z) \rightarrow 0 \qquad 
	\text{as } z \rightarrow \infty.
    \end{displaymath}
    Write 
    \( U^{n} = 
      \bigl( \begin{smallmatrix} 
         a_{n}&b_{n}\\ c_{n}&d_{n}
      \end{smallmatrix} \bigr)
    \),
    so
    \( U^{n}S = 
      \bigl( \begin{smallmatrix} 
        a_{n}&a_{n+1}\\ a_{n-1}&a_{n} 
      \end{smallmatrix} \bigr)
    \)
    and
    \begin{displaymath}
        \left( q \vert U^{n}S \right)(z)  
	    =  \left(a_{n-1}z-a_{n}\right)^{-2k} q\left(U^{n}Sz\right).
    \end{displaymath}
    Now \( \left(a_{n-1}z-a_{n}\right)^{-2k} \rightarrow 0 \) as
    \( z \rightarrow \infty \), since \( a_{n-1} \neq 0 \) for 
    \( 1 \leq n \leq p-2 \).
    We calculate that \( 0 < U^{n}S(\infty) < \infty \) 
    for \( 1 \leq n \leq p-2 \),
    and \( U^{n}S(\infty) \) is parabolic,
    so
    \( U^{n}S(\infty) \)
    cannot be a pole of \( q \), \emph{i.e.},
    \( q\left(U^{n}Sz\right) \) is bounded near \( \infty \).
    Thus 
    \( \left(a_{n-1}z-a_{n-2}\right)^{-2k} q\left(U^{n}Sz\right) 
        \rightarrow 0\) 
    as 
    \( z \rightarrow \infty \)
    \( (1 \leq n \leq p-2) \), and the claim holds.
    If \( q \) had a pole at \( \infty \) of order \( m>0 \),
    then 
    \[ q(z) = r(z) + c_{m}z^{m} + c_{m-1}z^{m-1} + \cdots + c_{0}, \]
    where
    \( r(\infty) = 0 \) and \( c_{m} \neq 0 \).
    Then
    \begin{displaymath}
        \left(q-q \vert S\right)(z) 
          = r(z) - r(z+\lambda) 
	     - m\lambda c_{m}z^{m-1} + \cdots ,
    \end{displaymath}
    where ``\(\cdots\)" denotes a polynomial
    of degree less than \( m-1 \).
    Thus as \( z \rightarrow \infty \),
    \begin{equation*}
	\left(q-q \vert S\right)(z) \rightarrow 
	\begin{cases}
	    -\lambda c_{1}, & 
	    \text{if } m=1, \\
	    \infty , & 
	    \text{if } m>1,
	\end{cases}
    \end{equation*}
    which contradicts the claim and establishes \emph{(i)}.

    Parts \emph{(ii)} and \emph{(iii)} follow immediately from 
    \emph{(i)} and the 
    first relation \eqref{eq:FirstRelation} in the form
    \begin{equation*}
        q\left(\frac{-1}{z}\right) = -z^{2k}q(z).
    \end{equation*}
\end{proof}

An RPF of weight \( 2k \) on \( G_{p} \) may have a pole \emph{only} at zero.
Meier and Rosenberger show in \cite{MR84} that
such an RPF is of the form
\begin{equation}
  q_{k,0}(z) = 
    \begin{cases}
      a_{0}(1-z^{-2k}), & \text{if } 2k \neq 2, \\
      a_{0}(1-z^{-2}) + b_{1}z^{-1}, & \text{if } 2k = 2.
    \end{cases}
    \label{eq:RPFPoleatZero}
\end{equation}

Given \( q \), an RPF of weight \( 2k \) on \( G_{p} \), we let
\(PP_{\alpha}[q] \) denote the principal part of \(q(z)\) at \(z=\alpha\).
Then by \eqref{eq:PoleDecomp2} and Lemma \ref{lem.PoleatZero}, 
\( q \) has the form
\begin{equation}
    q(z) = \sum_{\ell=1}^{L} 
           \sum_{\alpha \in Z_{\mathcal{A}_{\ell}} 
	         \cup TZ_{\mathcal{A}_{\ell}}}
               PP_{\alpha}[q](z) 
	       + b_{0}q_{k,0}(z)
	       + \sum_{n=1}^{2k-1}\frac{c_{n}}{z^{n}},
    \label{eq:RPFform1}
\end{equation}
where \( q_{k,0} \) is an RPF given by \eqref{eq:RPFPoleatZero}, 
\( b_{0} \) and each \( a_{n} \) and \( c_{n} \) is a constant,
and the \( \mathcal{A}_{\ell} \) are indexed hyperbolic \( G_{p} 
\)-equivalence classes of \( \lambda_{p} \)-BQFs.

\subsubsection{Nonzero Poles}
\label{subsubsec.nonzeropoles}

\begin{lemma}\label{lem.PoleOrder_k}
Let \(q\) be an RPF of weight \(2k\) on \(G_{p}\) with a nonzero pole at 
\(\alpha\).
Then the pole at \(\alpha\) is of order \(k\).
\end{lemma}

\begin{remark}
Schmidt \cite{Sch93} first states this Lemma with a note that the proof 
for RPFs on \(\Gamma(1)\) goes through.
The following proof generalizes the one for \(\Gamma(1)\)
by Choie and Zagier \cite{CZ93}.
\end{remark}

\begin{proof}
    Fix \( p \geq 3 \) and put \( U = U_{\lambda_{p}} \).
    Note that if \(f(z)\) has a pole at \(\beta\) and 
    \(V=\bigl( \begin{smallmatrix} *&*\\ c&d \end{smallmatrix} \bigr)\) 
    is a linear fractional transformation, then
    \(\left( f \mid V \right) (z) = (cz+d)^{-2k}f(Vz)\)
    has a pole at \(z=V^{-1}\beta\).
    In fact, 
    \begin{equation}
	PP_{V^{-1}\beta}\left[ f\mid V \right]
	= PP_{V^{-1}\beta}\left[ PP_{\beta}[f] \mid V \right].
	\label{eq:PPformulaV}
    \end{equation}
    By \eqref{eq:FirstRelation} 
    and \eqref{eq:PPformulaV} with \(V=T\),
    we have 
    \begin{equation}
	PP_{T\beta}[q] = -PP_{T\beta}\left[ PP_{\beta}[q] \mid T \right].
	\label{eq:PPformulaT}
    \end{equation}
    In a similar way, we use \eqref{eq:SecondRelation} 
    and \eqref{eq:PPformulaV} with \(V=U^{-t}\), 
    \( t \) an integer,
    to get 
    \begin{equation}
	PP_{U^{t}\beta}[q] 
	= -PP_{U^{t}\beta}\left[ PP_{\beta}[q] \mid U^{-t} \right].
	\label{eq:PPformulaUpower}
    \end{equation}

    Suppose that \( \alpha > 0 \).
    Then by the proof of Lemma \ref{lem.MinPoleSet},
    \(\alpha\) is fixed by an element of \( G_{p} \) of the form 
    \(M=U^{j_{r}}TU^{j_{r-1}}T \dots U^{j_{1}}T\),
    where \( 1 \leq j_{\nu} \leq p-1 \) for each \( \nu \geq 1 \).
    Applying \eqref{eq:PPformulaT} and \eqref{eq:PPformulaUpower} \(r\) 
    times each, we get
    \begin{equation*}
	PP_{M\alpha}[q] 
	= PP_{M\alpha}\left[ PP_{M\alpha}[q] \mid M^{-1} \right],
    \end{equation*}
    and, since \(M\alpha=\alpha\),
    \begin{equation}
	PP_{\alpha}[q] 
	= PP_{\alpha}\left[ PP_{\alpha}[q] \mid M^{-1} \right].
	\label{eq:PP_slash_eqn}
    \end{equation}
    If \( \alpha < 0 \), then \( T\alpha \in \Pset^{+} \) is fixed by 
    an element of \( G_{p} \) of the form
    \( N = U^{j_{r}}TU^{j_{r-1}}T \dots U^{j_{1}}T \).
    Thus \( \alpha \) is fixed by 
    \( M = TNT = TU^{j_{r}}TU^{j_{r-1}}T \dots U^{j_{1}} \).
    We apply \eqref{eq:PPformulaT} and \eqref{eq:PPformulaUpower} \(r\) 
    times each in this case as well, 
    and \( q \) also satisfies 
    \eqref{eq:PP_slash_eqn} when \( \alpha < 0 \).

    Let \(m\) be the order of the pole of \(q\) at \(\alpha\) 
    and suppose that \(q\) is normalized so that 
    \(PP_{\alpha}[q](z)=(z-\alpha)^{-m} + r_{\alpha}(z) \),
    where \( r_{\alpha}(z) \) is a rational function 
    with a pole at \(\alpha\) of order less than \(m\).
    Put \(M=\bigl( \begin{smallmatrix} a&b\\ c&d \end{smallmatrix} \bigr)\)
    so \(M^{-1}=\bigl( \begin{smallmatrix} d&-b\\ -c&a \end{smallmatrix} 
    \bigr)\).
    We calculate
    \begin{eqnarray*} 
	\left( PP_{\alpha}[q] \mid M^{-1} \right) (z) 
	& = & (-cz+a)^{-2k} \left(\left( M^{-1}z-\alpha \right)^{-m} 
	        + r_{\alpha}\left( M^{-1}z \right)\right)          \\
	& = & (-cz+a)^{m-2k} \left( (c\alpha + d)z - (a\alpha + b) \right)^{-m} 
	    + r_{\alpha,\frac{a}{c}}(z) \\
	& = & \frac{(-cz+a)^{m-2k}}{(c\alpha + d)^{m}} 
	    \left(z - M\alpha \right)^{-m} + r_{\alpha,\frac{a}{c}}(z) \\
	& = & \frac{(-c\alpha +a)^{m-2k}}{(c\alpha + d)^{m}} (z - \alpha)^{-m} 
	    + \tilde{r}_{\alpha,\frac{a}{c}}(z) \\
	& = & (c\alpha + d)^{2k-2m} (z - \alpha)^{-m} 
	    + \tilde{r}_{\alpha,\frac{a}{c}}(z),
    \end{eqnarray*} 
    where \( r_{\alpha,\frac{a}{c}} \) and 
    \( \tilde{r}_{\alpha,\frac{a}{c}} \) are 
    rational functions, each with a pole at \( \alpha \) of order less than \( m \)
    and a pole at \( \frac{a}{c} \).
    We have used partial fractions, the fact that \(M\) fixes \(\alpha\),
    and the fact that \((-c\alpha + a)=(c\alpha + d)^{-1}\).
    Then by \eqref{eq:PP_slash_eqn} we have  
    \begin{eqnarray*} 
	PP_{\alpha}[q](z)
	& = & PP_{\alpha} \left[ (c\alpha + d)^{2k-2m} (z - \alpha)^{-m} 
	    + \tilde{r}_{\alpha,\frac{a}{c}}(z) \right]     \\
	& = & (c\alpha + d)^{2k-2m} (z - \alpha)^{-m} 
	    + \tilde{r}_{\alpha}(z),
    \end{eqnarray*} 
    where \( \tilde{r}_{\alpha} \) is a rational function with a pole at 
    \( \alpha \) of order less than \( m \).
    Thus \((c\alpha + d)^{2k-2m} = 1\),
    so either \(m=k\) or \(c\alpha +d=1\).
    But if \(c\alpha +d=1\), then \(\alpha = \frac{1-d}{c}\)
    is parabolic, a contradiction.
    We conclude that \(m=k\).
\end{proof}

\begin{lemma}\label{lem.UniquePP}
Suppose that \(\alpha \neq 0\) can occur as a pole of an RPF of weight 
\(2k\) on \(G_{p}\).  
Then there exists a unique function of the form 
\begin{equation*}
    q_{k,\alpha}(z) = \frac{1}{(z-\alpha)^{k}} 
			+ \frac{a_{1}(k,\alpha)}{(z-\alpha)^{k-1}}
                  + \dots + \frac{a_{k-1}(k,\alpha)}{z-\alpha},
\end{equation*}
such that for \emph{any} 
RPF \(q\) of weight \(2k\) on \(G_{p}\),
\(PP_{\alpha}[q]\) is a constant multiple of \(q_{k,\alpha}\).
\end{lemma}

\begin{remark}\label{rem:PP[q]}
    The function \(q_{k,\alpha}\) is uniquely determined by 
    \(k\) and \(\alpha\), 
    but is independent of the rational period function \(q\).
\end{remark}

\begin{proof}
Functions of the given form exist,
since \(PP_{\alpha}[q]\) has a pole of order 
\( k \) at \( \alpha \) and can be normalized so 
that the coefficient of \( (z-\alpha)^{-k} \) is \( 1 \).
For uniqueness, suppose that \(q_{k,\alpha}\) and \(r_{k,\alpha}\) 
are functions satisfying the hypotheses of the Lemma, 
with \(q_{k,\alpha} \neq r_{k,\alpha}\).
Then there exist RPFs \( q \) and \( r \)
and nonzero constants \( a \) and \( b \) such that
\(PP_{\alpha}[q] = aq_{k,\alpha}\) and 
\(PP_{\alpha}[r] = br_{k,\alpha}\).
But then \( q-\frac{a}{b}r \) is an RPF with a pole at
\( \alpha \) of order less than \( 2k \),
a contradiction.
\end{proof}

The following Lemma gives a formula for \( q_{k,\alpha} \).
\begin{lemma}\label{lem.ExplicitPP}
    Suppose that \( \alpha \neq 0 \) can occur as a pole of an RPF 
    \( q \) of weight \( 2k \) on \( G_{p} \).
    Then
    \begin{equation}
	q_{k,\alpha}(z)
	= PP_{\alpha} 
	  \left[ \frac{D^{k/2}}{Q_{\alpha}(z,1)^{k}} \right]
	= PP_{\alpha} 
	         \left[ \frac{(\alpha-\alpha^{\prime})^{k}}                                       
		 {(z-\alpha)^{k}(z-\alpha^{\prime})^{k}}
                 \right],
    \label{eq:PPatalpha}
    \end{equation}
    where \(Q_{\alpha}\) is the \(\lambda_{p}\)-BQF associated to 
    \(\alpha\), 
    \(D\) is the discriminant of \( Q_{\alpha} \),
    and \( \alpha^{\prime} \) is the Hecke \( \lambda_{p} \)-conjugate 
    of \( \alpha \).
\end{lemma}

\begin{remark}
Schmidt \cite{Sch93} states this Lemma in a modified form 
for \( \alpha > 0 \) and asserts that 
the proof follows from the proof in the classical case.
The following proof holds for
\( \alpha \) positive or negative.
\end{remark}

\begin{proof}
    Write \(Q_{\alpha}(z) = Q_{\alpha}(z,1) = Az^{2}+Bz+C \).
    Let 
    \(M=\bigl( \begin{smallmatrix} a&b\\ c&d \end{smallmatrix} \bigr) \) 
    be an element of \(G_{p}\) which fixes both \(\alpha\) and 
    \(\alpha^{\prime}\).
    We calculate
    \begin{eqnarray*} 
	\lefteqn{ \left( Q_{\alpha}^{-k} \mid M^{-1} \right) (z) }  \\
	& = & 
	(-cz+a)^{-2k} A^{-k} \left( \frac{dz-b}{-cz+a}-\alpha \right)^{-k}
	\left( \frac{dz-b}{-cz+a}-\alpha^{\prime} \right)^{-k} \\
	& = & 
	A^{-k} \left( (c\alpha + d)z - (a\alpha + b) \right)^{-k} 
	\left( (c\alpha^{\prime} + d)z - (a\alpha^{\prime} + b) \right)^{-k} \\
	& = & 
	A^{-k} (c\alpha + d)^{-k}(c\alpha^{\prime} + d)^{-k}
	\left(z - \frac{a\alpha + b}{c\alpha + d} \right)^{-k}
	\left(z - \frac{a\alpha^{\prime} + b}{c\alpha^{\prime} + d} \right)^{-k} 
	\\
	& = & A^{-k}(z-\alpha)^{-k}(z-\alpha^{\prime})^{-k}   \\
	& = & Q_{\alpha}(z)^{-k}.
    \end{eqnarray*} 
    We have used the fact that \((c\alpha + d)(c\alpha^{\prime} + d)=1\).
    From this it follows that
    \( Q_{\alpha}^{-k} \) satisfies \eqref{eq:PP_slash_eqn}, 
    as does \( q \).
    By an argument similar to the proof of Lemma \ref{lem.UniquePP},
    \(PP_{\alpha}[Q_{\alpha}^{-k}]\) is a nonzero multiple of \(PP_{\alpha}[q]\),
    which is in turn a multiple of \(q_{k,\alpha}\).
    Thus
    \(q_{k,\alpha}\) is a multiple of \(PP_{\alpha}[Q_{\alpha}^{-k}]\).
    Finally 
    \begin{eqnarray*}
	PP_{\alpha}[Q_{\alpha}^{-k}](z)
	& = & PP_{\alpha}[A^{-k}(z-\alpha)^{-k}(z-\alpha^{\prime})^{-k}]   \\
	& = & A^{-k}(z-\alpha)^{-k}(\alpha-\alpha^{\prime})^{-k} + \dotsb   \\
	& = & D^{-k/2}(z-\alpha)^{-k} + \dotsb,
    \end{eqnarray*}
    so the proportionality constants are as given in the Lemma.
\end{proof}

We now give a complete description of any RPF of weight \( 2k \) on 
\( G_{p} \).

\begin{theorem}\label{thm.AnyRPFdiscription}
    An RPF of weight \( 2k \in 2\Z \) on \( G_{p} \) is of the form
    \begin{equation*}
        q(z) = \sum_{\ell=1}^{L} C_{\ell} 
    	    \sum_{\alpha \in Z_{\mathcal{A}_{\ell}}}
	    \left( q_{k,\alpha}(z) - q_{k,T\alpha}(z) \right) 
	    + c_{0}q_{k,0}(z)
	    + \sum_{n=1}^{2k-1}\frac{c_{n}}{z^{n}},
    \end{equation*}
    where each
    \( \mathcal{A}_{\ell} \) is a \( G_{p} \)-equivalence class of 
        \( \Z[\lambda_{p}] \)-BQFs,
    \( Z_{\mathcal{A}_{\ell}} \) is the cycle of positive poles 
        associated with \( \mathcal{A}_{\ell} \),
    \( q_{k,\alpha} \) is given by \eqref{eq:PPatalpha},
    \( q_{k,0} \) is given by \eqref{eq:RPFPoleatZero},
    and the \( C_{\ell} \) and \( c_{n} \) are all constants.
\end{theorem}

\begin{proof}
    By \eqref{eq:RPFform1} and Lemma \ref{lem.UniquePP} we have that
    any RPF \( q \) of weight \(2k\) on \(G_{p}\) has the form
    \begin{equation*}
	q(z) = \sum_{\ell=1}^{L}
	    \sum_{\alpha \in Z_{\mathcal{A}_{\ell}} 
                \cup TZ_{\mathcal{A}_{\ell}} } 
	        C_{\alpha}q_{k,\alpha}(z) 
	    + c_{0}q_{k,0}(z)
	    + \sum_{n=1}^{2k-1}\frac{c_{n}}{z^{n}},
    \end{equation*}
    where each \( \mathcal{A}_{\ell} \) is a hyperbolic 
    \( G_{p} \)-equivalence class of \( \lambda_{p} \)-BQFs and 
    \( q_{k,0} \) is given by \eqref{eq:RPFPoleatZero}.
    From \eqref{eq:PPformulaT} and \eqref{eq:PPformulaUpower}
    we see that the coefficients \(C_{\alpha}\)
    alternate in sign as \( \alpha \) goes through any cycle
    \(\left\{ \alpha_{1}, T\alpha_{1}, U^{j_{1}}T\alpha_{1}, 
    TU^{j_{1}}T\alpha_{1}, 
    \dots \right\} \).
    Thus for each irreducible system of poles 
    \(Z_{\mathcal{A}_{\ell}} \cup TZ_{\mathcal{A}_{\ell}}\)
    there is a constant \(C_{\ell}\)
    such that \(C_{\alpha}=C_{\ell}\) 
    for each \(\alpha \in Z_{\mathcal{A}_{\ell}}\)
    and \(C_{\alpha}=-C_{\ell}\) 
    for each \(\alpha \in TZ_{\mathcal{A}_{\ell}}\).
    As a result any RPF of weight \(2k\) on \(G_{p}\) has the 
    given form.
\end{proof}

Lemma \ref{lem.NegPoles} allows us to identify negative poles using Hecke conjugation instead of the 
action of \( T \) in the Corollary.
\begin{corollary}
\label{cor:AnyRPFdiscription}
    An RPF of weight \( 2k \in 2\Z \) on \( G_{p} \) is of the form
    \begin{equation}
	\label{eq:RPFform3}
        q(z) = \sum_{\ell=1}^{L} C_{\ell} 
    	    \left(\sum_{\alpha \in Z_{\mathcal{A}_{\ell}}}
	            q_{k,\alpha}(z) 
		- \sum_{\alpha \in Z_{-\mathcal{A}_{\ell}}}
		    q_{k,\alpha^{\prime}}(z)
	    \right) 
	    + c_{0}q_{k,0}(z)
	    + \sum_{n=1}^{2k-1}\frac{c_{n}}{z^{n}},
    \end{equation}
    where each
    \( \mathcal{A}_{\ell} \) is a \( G_{p} \)-equivalence class of 
        \( \Z[\lambda_{p}] \)-BQFs,
    \( Z_{\mathcal{A}_{\ell}} \) is the cycle of positive poles 
        associated with \( \mathcal{A}_{\ell} \),
    \( q_{k,\alpha} \) is given by \eqref{eq:PPatalpha},
    \( q_{k,0} \) is given by \eqref{eq:RPFPoleatZero},
    and the \( C_{\ell} \) and \( c_{n} \) are all constants.
\end{corollary}

\section{HECKE-SYMMETRIC POLE SETS, $k$ ODD}
\label{sec.symmetric,kodd}

The following theorem characterizes of RPFs 
of weight \( 2k \) on \( G_{p} \), for \( k \) odd, with
Hecke-symmetric sets of poles.
\begin{theorem}\label{thm.SymmetricCharacterization}
    Suppose that \( k \) is odd.
    
    Suppose \( q \) is an RPF of weight \( 2k \) on \( G_{p} \)
    with Hecke \( \lambda_{p} \)-symmetric
    irreducible systems of poles.
    Then
    \( q \) is of the form
    \begin{equation}
        q(z) = \sum_{\ell=1}^{L} d_{\ell}
               \sum_{\alpha \in Z_{\mathcal{A}_{\ell}}} 
	       Q_{\alpha}(z,1)^{-k}
	     + c_{0}q_{k,0}(z),
        \label{eq:RPFsymmetrickodd}
    \end{equation}
    where each
    \( \mathcal{A}_{\ell} \) 
    is a \( G_{p} \)-equivalence class of 
    \( \Z[\lambda_{p}] \)-BQFs
    satisfying \( -\mathcal{A}_{\ell} = \mathcal{A}_{\ell} \),
    the \( d_{\ell} \) \( (1 \leq \ell \leq M) \) are constants,
    and \( q_{k,0}(z) \) is given by \eqref{eq:RPFPoleatZero}.
    
    Conversely,
    any rational function of the form
    \eqref{eq:RPFsymmetrickodd}
    is an RPF of weight \( 2k \) 
    on \( G_{p} \) with 
    Hecke \( \lambda_{p} \)-symmetric irreducible systems of poles.
\end{theorem}

\begin{proof}
    Suppose that \( q \) is an RPF of weight \( 2k \) on \( G_{p} \)
    with Hecke \( \lambda_{p} \)-symmetric
    irreducible systems of poles.
    Then \( -\mathcal{A} = \mathcal{A} \)
    for each \( G_{p} \)-equivalence class of \( \lambda_{p} \)-BQFs 
    associated with poles of \( q \).
    Thus the expression for \( q \) in \eqref{eq:RPFform3} simplifies 
    to
    \begin{equation}
        q(z) = \sum_{\ell=1}^{L} C_{\ell}
                   \sum_{\alpha \in Z_{\mathcal{A}_{\ell}}} 
	           \left(q_{k,\alpha}(z)-q_{k,\alpha^{\prime}}(z)\right)
	       + c_{0}q_{k,0}(z)
	       + \sum_{n=1}^{2k-1} \frac{c_{n}}{z^{n}},
        \nonumber
    \end{equation}
    where 
    \( q_{k,0}(z) \) is given by \eqref{eq:RPFPoleatZero} and 
    \( C_{\ell} \) \( (1 \leq \ell \leq L) \)
    and 
    \( c_{n} \) \( (1 \leq n \leq 2k-1) \) are all constants.
    A calculation shows that for any pole \( \alpha \),
    \begin{equation}
        q_{k,\alpha}(z)-q_{k,\alpha^{\prime}}(z)
	= PP_{\alpha} 
	  \left[ \frac{D^{k/2}}{Q_{\alpha}(z,1)^{k}} \right]
	  - (-1)^{k}
	  PP_{\alpha^{\prime}} 
	  \left[ \frac{D^{k/2}}{Q_{\alpha}(z,1)^{k}} \right],
	\nonumber
    \end{equation}
    where \( D \) is the discriminant of \( Q_{\alpha} \).
    Since \( k \) is odd,
    \begin{equation}
        q_{k,\alpha}(z)-q_{k,\alpha^{\prime}}(z)
	= \frac{D^{k/2}}{Q_{\alpha}(z,1)^{k}},
	\nonumber
    \end{equation}
    so  
    \begin{equation}
        q(z) = \sum_{\ell=1}^{L} d_{\ell}
                   \sum_{\alpha \in Z_{\mathcal{A}_{\ell}}} 
	           Q_{\alpha}(z,1)^{-k}
	       + c_{0}q_{k,0}(z)
	       + \sum_{n=1}^{2k-1} \frac{c_{n}}{z^{n}},
	\label{eq:RPFform4}
    \end{equation}
    where 
    \( d_{\ell} = C_{\ell}D_{\ell}^{k/2} \) and
    \( D_{\ell} \) is the discriminant of 
    \( \lambda_{p} \)-BQFs in \( \mathcal{A}_{\ell} \).

    Schmidt proves in \cite[Theorem 3.1]{Sch93} that if 
    \( k \) is odd, then 
    \begin{equation}
        r(z) = \sum_{\alpha \in Z_{\mathcal{A}} 
	             \cup Z_{-\mathcal{A}}}
	       Q_{\alpha}(z,1)^{-k}
	\nonumber
    \end{equation}
    is an RPF of weight \( 2k \) on \( G_{p} \).
    Thus the first sum in \eqref{eq:RPFform4} is an
    RPF of weight \( 2k \) on \( G_{p} \).
    Since \( q \) and \( q_{k,0} \) are also RPFs,
    we must have that 
    \( \sum_{n=1}^{2k-1} \frac{c_{n}}{z^{n}} \)
    is an RPF of weight \( 2k \) on \( G_{p} \).
    But then we must have 
    \( c_{n}=0 \) for \( 1 \leq n \leq 2k-1 \)
    or else we contradict Lemma \ref{lem.PoleatZero} (\emph{ii}).
    Thus \( q \) has the form \eqref{eq:RPFsymmetrickodd}.

    For the converse,
    assume that \( q \) has the form 
    \eqref{eq:RPFsymmetrickodd}.
    Since 
    \( \sum_{\ell=1}^{L} d_{\ell}
               \sum_{\alpha \in Z_{\mathcal{A}_{\ell}}} 
	       Q_{\alpha}(z,1)^{-k} \)
    and \( q_{k,0} \) are both 
    RPFs of weight \( 2k \) on \( G_{p} \),
    \( q \) is also an RPF of weight \( 2k \) on \( G_{p} \).
    The set of nonzero poles of \( q \) is
    \begin{equation}
        P^{*}(q) = \bigcup_{\ell=1}^{L}
                       \left( Z_{\mathcal{A}_{\ell}} \cup 
                       Z_{\mathcal{A}_{\ell}}^{\prime} \right),
        \nonumber
    \end{equation}
    which is clearly Hecke-symmetric.

\end{proof}

\end{document}